# Lie Symmetry-Guided Adaptive Sampling for Physics-Informed Neural Networks

RUNZE ZHAO

## Abstract

Physics-informed neural networks (PINNs) have emerged as a promising mesh-free framework for solving partial differential equations (PDEs), with broad applications in physics, engineering, and medicine. However, their performance is highly sensitive to the sampling strategy for training points. Conventional uniform sampling often fails to adequately resolve regions with complex solution features, leading to low efficiency and limited accuracy. To address these challenges, we propose a Lie symmetry-driven adaptive sampling framework for PINNs, termed Lie-AS-ISC. For PDEs admitting Lie symmetry groups, symmetry transformations are exploited to derive intrinsic functional relationships between the solution and its variables, enabling the construction of basis points with exact solutions that span the solution domain. Building on this, an adaptive sampling strategy is developed by jointly considering absolute error and invariant surface conditions, while a dual-index loss is introduced to enhance training constraints. The proposed method significantly improves prediction accuracy. The experimental results show that this method consistently outperforms existing methods in accuracy.

## 1. Introduction

Partial differential equations (PDEs) constitute a fundamental class of mathematical models and are widely used in fields such as medicine [1], weather forecasting [2], and fluid dynamics [3]. However, for most PDEs, deriving analytical solutions is intractable, making their efficient and accurate solution a longstanding challenge. Traditional approaches rely on numerical methods, including the finite difference method [4], finite element method [5], and finite volume method [6], which have been extensively employed to approximate solutions when analytical forms are unavailable. Despite their success, these methods typically require the construction of computational grids, which becomes increasingly expensive and inefficient when dealing with high-dimensional problems.

Deep neural networks have achieved remarkable success across a wide range of scientific and engineering domains, including medicine [7], biology [8], object recognition [9], and disease discovery [10]-[12]. In recent years, deep learning has also shown great potential in solving differential equations, leading to the development of several novel approaches such as the deep Galerkin method (DGM) [13] and the Deep Ritz method [14]. Among these, physics-informed neural networks (PINNs) [15] have emerged as a representative framework. As a mesh-free approach, PINNs can significantly reduce computational costs and offer promising application potential. However, their performance is often hindered by issues such as training instability, limited accuracy, and low efficiency, leaving substantial room for further improvement.

In recent years, numerous efforts have been devoted to improving the performance of PINNs. Jagtap et al. [16] introduced the extended physics-informed neural network (XPINN), which decomposes the solution domain into multiple non-overlapping subdomains and trains separate networks on each subdomain, thereby enhancing representation capability and enabling parallel computation. Hou et al. [17] proposed AMAW-PINN, which integrates adaptive sampling with adaptive loss weighting by using residuals and gradients as indicators, allowing denser sampling in regions that are difficult to approximate. Inspired by the phototactic behavior of fireflies, Wang et al. [18] developed FAMAW, a bio-inspired approach that guides training points toward regions with large residuals or gradients. Cao et al. [19] proposed ARPINN, an adaptive residual-splitting method that decomposes the PDE residual into multiple components and dynamically adjusts their weights to improve training efficiency. In addition to these representative approaches, a variety of other strategies have been developed to address different limitations of PINNs [20]-[32].

Despite these advances, several limitations remain. Existing adaptive sampling strategies that rely on residuals or gradients as indicators can partially alleviate poor performance in challenging regions; however, such indirect criteria are often insufficient to fully capture complex solution behaviors.

To address these issues, we propose a Lie symmetry-based adaptive sampling framework for PINNs, termed Lie-AS-ISC. Specifically, symmetry-induced points derived from the Lie symmetries of the governing equations are employed as global reference points. Based on these, adaptive sampling is performed by jointly considering absolute error and invariant surface conditions, and a dual-index loss is imposed to further strengthen training constraints. By directly incorporating intrinsic equation symmetries into the sampling process, the proposed method provides more effective guidance in challenging regions and significantly improves model performance. Numerical experiments demonstrate notable gains in both accuracy and efficiency. The remainder of this paper is organized as follows. Section 2 introduces the relevant background, Section 3 presents the proposed method, Section 4 reports numerical experiments, and Section 5 concludes the paper.

## 2. Background knowledge

### 2.1 PINN

For a given PDE:

$$h(x,u(x))=0, \tag{1}$$

Where

$$h(x,u(x))=h\left(x;\frac{\partial u}{\partial x_1},\cdots,\frac{\partial u}{\partial x_n};\frac{\partial^2 u}{\partial x_1\partial x_2},\cdots,\frac{\partial^2 u}{\partial x_1\partial x_n},\cdots,\lambda\right), \tag{2}$$

$\mathrm{B}(u(x),x)=0$ denotes the boundary condition, and $u(x)$ is the solution at $x$. In forward problems, the parameter $\lambda$ is known while the solution $u(x)$ is unknown; in inverse problems, the parameter $\lambda$ is inferred from data. PINN can address both forward and inverse problems. For the forward problem, a neural network $\hat{u}(x,\theta)$ represents solution $u$, and training parameter $\theta$ is used to approximate solution $u$, The network is trained by minimizing the following loss function::

$$L_h(\theta;N_h)=\frac{1}{|N_h|}\sum_{x\in N_h}\left|h\left(x;\frac{\partial \hat{u}}{\partial x_1},\cdots,\frac{\partial \hat{u}}{\partial x_n};\frac{\partial^2 \hat{u}}{\partial x_1\partial x_2},\cdots,\frac{\partial^2 \hat{u}}{\partial x_1\partial x_n},\cdots,\lambda\right)\right|^2, \tag{3}$$

$$L_b\left(\theta;N_b\right)=\frac{1}{|N_b|}\sum_{x\in N_b}\left|B\left(\hat{u},x\right)\right|^2, \tag{4}$$

$$L\left(\theta;N\right)=\omega_h L_h\left(\theta;N_h\right)+\omega_b L_b\left(\theta;N_b\right), \tag{5}$$

Here, $L_h\left(\theta;N_h\right)$ and $L_b\left(\theta;N_b\right)$ denote the residual loss and boundary condition loss, respectively, while $\omega_h$ and $\omega_b$ are weighting coefficients that balance the two loss terms. $N_h$ denotes the number of collocation points in the interior domain, and $N_b$ denotes the number of collocation points on the boundary[15].

For inverse problems, an additional observation data loss term is introduced:

$$L\left(\theta,\lambda;N_i\right)=\frac{1}{|N_i|}\sum_{x\in N_i}\left|u\left(\hat{x}\right)-u\left(x\right)\right|^2, \tag{6}$$

Accordingly, the total loss function is defined as

$$L\left(\theta,\lambda;N\right)=\omega_h L_h\left(\theta,\lambda;N_h\right)+\omega_b L_b\left(\theta,\lambda;N_b\right)+\omega_i L_i\left(\theta,\lambda;N_i\right), \tag{7}$$

Where $L_i$ denotes the data loss term and $\omega_i$ is its corresponding weight.

## 2.2 Lie symmetric group

Many partial differential equations possess symmetries, and symmetry analysis is an important method for investigating the properties of partial differential equations. For a one-parameter Lie group, the corresponding infinitesimal transformations can be expressed as

$$x_{new}=x+a\delta+O\left(\varepsilon^2\right),t_{new}=t+b\delta+O\left(\varepsilon^2\right),u_{new}=u+c\delta+O\left(\varepsilon^2\right), \tag{8}$$

If a PDE admits a Lie symmetry group, its form remains invariant under the corresponding symmetry transformations [33]-[34].

Consider a PDE of the form:

$$f:=u_t+\mathcal{N}\left[u\right]=0, \tag{9}$$

A Lie symmetry group associated with this equation can be characterized by an infinitesimal generator

$$\mathcal{X} = a(x,t,u)\partial_x + b(x,t,u)\partial_t + c(x,t,u)\partial_u, \tag{10}$$

where $a = a(x,t,u), b = b(x,t,u), c = c(x,t,u)$ are smooth functions. The operator $\mathcal{X}$ represents the infinitesimal generator of the symmetry group [33].

The infinitesimal invariance criterion for the equation is given by:

$$pr^{(k)}\mathcal{X}\left(u_t + \mathcal{N}[u]\right) = 0, \tag{11}$$

where $pr^{(k)}\mathcal{X}$ denotes the $k-th$ prolongation of the infinitesimal generator $\mathcal{X}$, defined as

$$pr^{(k)}\mathcal{X} = \mathcal{X} + \eta_t^{(1)}\partial_{u_t} + \eta_x^{(1)}\partial_{u_1} + \eta_x^{(2)}\partial_{u_2} + \cdots + \eta_x^{(k)}\partial_{u_k}. \tag{12}$$

Here $u_k = \dfrac{\partial^k u}{\partial x^k}, k = 1,2,\cdots,$ and the coefficients $\eta_t^{(1)}$ and $\eta_x^{(k)}$ are determined as follows:

$$\begin{aligned} \eta_t^{(1)} &= D_t\left(\eta - bu_t - au_x\right) + bD_t^2 u + aD_x D_t u, \\ \eta_x^{(k)} &= D_x^k\left(\eta - bu_t - au_x\right) + bD_t D_x^k u + aD_x^{k+1} u, \end{aligned} \tag{13}$$

Here, $D_t$ and $D_x$ denote the total derivative operators with respect to $t$ and $x$, respectively, defined as[35]

$$\begin{aligned} D_t &= \partial_t + u_t\partial_u + u_{xt}\partial_{u_x} + u_{tt}\partial_{u_t} + \cdots, \\ D_x &= \partial_x + u_1\partial_u + u_2\partial_u + u_{xt}\partial_{u_t} + \cdots, \end{aligned} \tag{14}$$

By solving for the infinitesimal coefficients $a$, $b$, and $c$, the infinitesimal generator $\mathcal{X}$ can be determined, from which the corresponding Lie symmetry group is derived. Furthermore, the invariant surface conditions can be obtained based on the symmetry group [33].

# 3. Methods

## 3.1 Base points with exact solutions

Inspired by the work of Liu et al. on domain decomposition[35], we extend the interfaces generated via symmetric transformations to the entire domain to acquire global base points.

For PDEs admitting a translational symmetry group, the one-parameter transformation can be expressed as

$$G_1 : x^* = x + \xi\varepsilon, t^* = t + \eta\varepsilon, u^* = u + \tau\varepsilon, \tag{15}$$

The corresponding invariant surface condition is given by:

$$g := \tau - \xi u_x - \eta u_t = 0, \tag{16}$$

From the transformation relations:

$$\begin{cases} x^* = x + \xi\varepsilon \\ t^* = t + \eta\varepsilon \end{cases}, \tag{17}$$

The relationship between $t^*$ and $x^*$ can be derived as:

$$x^* = \frac{\xi}{\eta} t^* + \frac{\eta x - \xi t}{\eta}, \eta \neq 0 \tag{18}$$

From the transformation relations:

$$\begin{cases} x^* = x + \xi\varepsilon \\ u^* = u + \tau\varepsilon \end{cases}, \tag{19}$$

The relationship between $x^*$ and $u^*$ can be derived as:

$$u^* = \frac{\tau}{\xi} x^* + \frac{\xi u - \tau x}{\xi}, \xi \neq 0 \tag{20}$$

From Eqs. (18) and (20), it follows that, given a point $A(x,t,u)$ with an exact solution value, a functional relationship among the transformed variables $x^*, t^*, u^*$ can be established. In particular, when either $\xi = 0$ or $\eta = 0$, the relation between $x^*$ and $t^*$ reduces to a linear form. Since points on the initial and boundary conditions possess exact solution values, each such point induces a corresponding relationship among $x^*, t^*, u^*$. By controlling the number of these points, the number of induced relationships can be effectively regulated. When a sufficient number of relationships is obtained, their images collectively cover a large portion of the solution domain. On these invariant mappings, a large number of points with exact solution values can be sampled, either uniformly or non-uniformly, thereby yielding base points that span the entire domain. Figure 1 illustrates the distribution of base points in the solution domain for $\xi = 3, \eta = 1$ and $\tau = 0$, with $x \in [0,1]$ and $t \in [0,1]$, under different numbers of sampling points on the initial and boundary conditions.

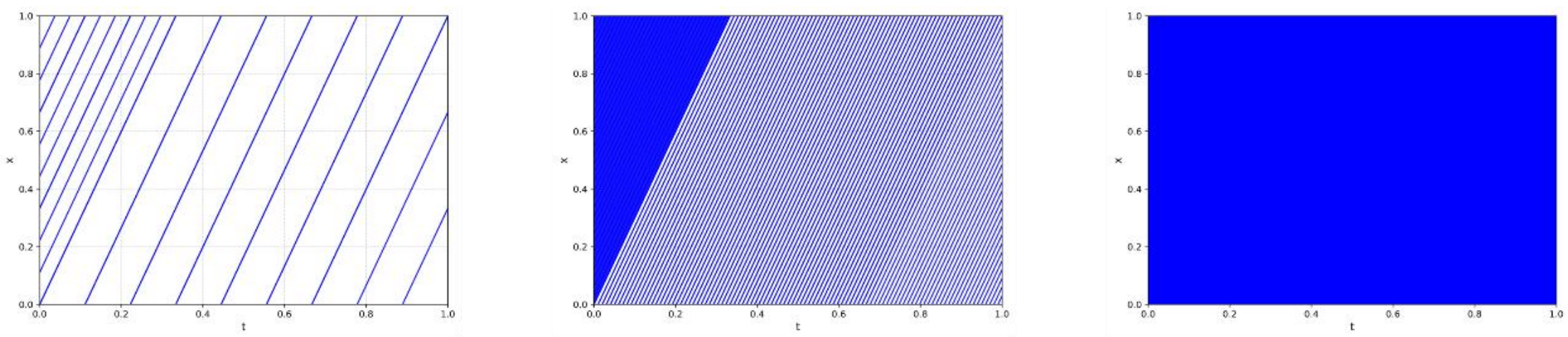

N=10 N=100 N=1000

Fig. 1. Distribution of base points in the solution domain obtained by sampling 10, 100, and 1000 points on the initial and boundary conditions, respectively.

It can be observed that as the number of sampled points increases, the base points progressively cover the solution domain. When 1000 points are sampled on the initial and boundary conditions, the base points nearly span the entire domain.

Next, consider PDEs admitting a rotational symmetry group. The corresponding one-parameter transformation is given by

$$G_2 : x^* = x\cos\varepsilon - t\sin\varepsilon, t^* = x\sin\varepsilon + t\cos\varepsilon, u^* = u, \tag{21}$$

The associated invariant surface condition is

$$tu_x - xu_t = 0, \tag{22}$$

From the transformation relations

$$\begin{cases} x^* = x\cos\varepsilon - t\sin\varepsilon \\ t^* = x\sin\varepsilon + t\cos\varepsilon \\ u^* = u \end{cases}, \tag{23}$$

The relationship among $x^*, t^*$ and $u^*$ can be derived accordingly:

$$\begin{cases} x^{*2} + t^{*2} = x^2 + t^2 \\ u^* = u \end{cases}, \tag{24}$$

Following the same procedure, the distribution of base points in the solution domain can be obtained for $x \in [0,1]$ and $t \in [0,1]$, as illustrated in Fig. 2.

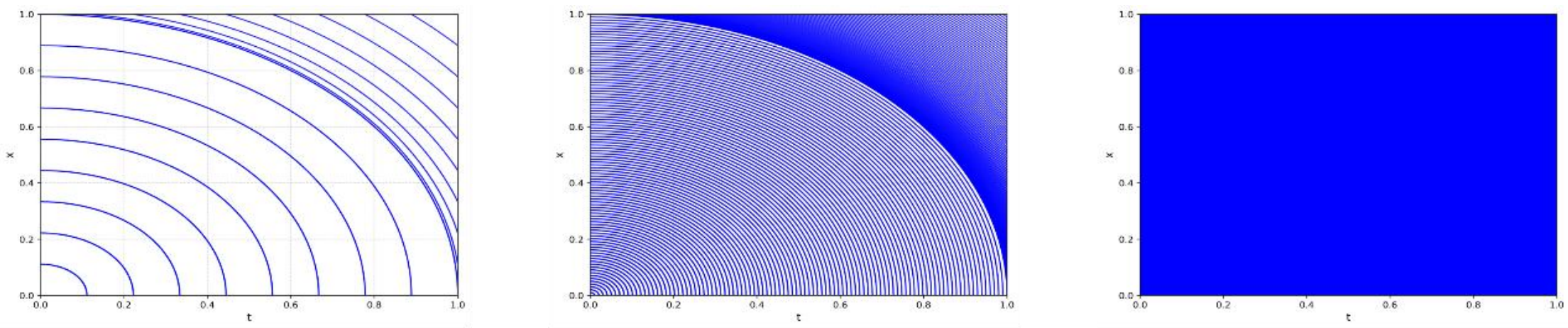

N=10 N=100 N=1000

Fig. 2. Distribution of base points in the solution domain obtained by sampling 10, 100, and 1000 points on the initial and boundary conditions, respectively.

Finally, PDEs admitting a scaling symmetry group. The corresponding transformation is given by

$$G_3 : x^* = \lambda^j x, t^* = \lambda^k t, u^* = \lambda^l u, \tag{25}$$

The associated invariant surface condition is

$$jxu_x + ktu_t = lu, \tag{26}$$

From the transformation relations

$$\begin{cases} x^* = \lambda^j x \\ t^* = \lambda^k t \end{cases}, \begin{cases} x^* = \lambda^j x \\ u^* = \lambda^l u \end{cases}, \tag{27}$$

The relationship among $x^*, t^*$ and $u^*$ can be derived accordingly:

$$\begin{cases} x^* = \dfrac{x}{t^{\frac{j}{k}}} t^{*\frac{j}{k}} \\ u^* = \dfrac{u}{x^{\frac{l}{j}}} x^{*\frac{l}{j}} \end{cases}, \tag{28}$$

Following the same procedure, the distribution of base points in the solution domain can be obtained. As an illustrative example, we consider $j = 2, k = 1, l = 0$, with $x \in [0.5,1]$ and $t \in [0.5,1]$.

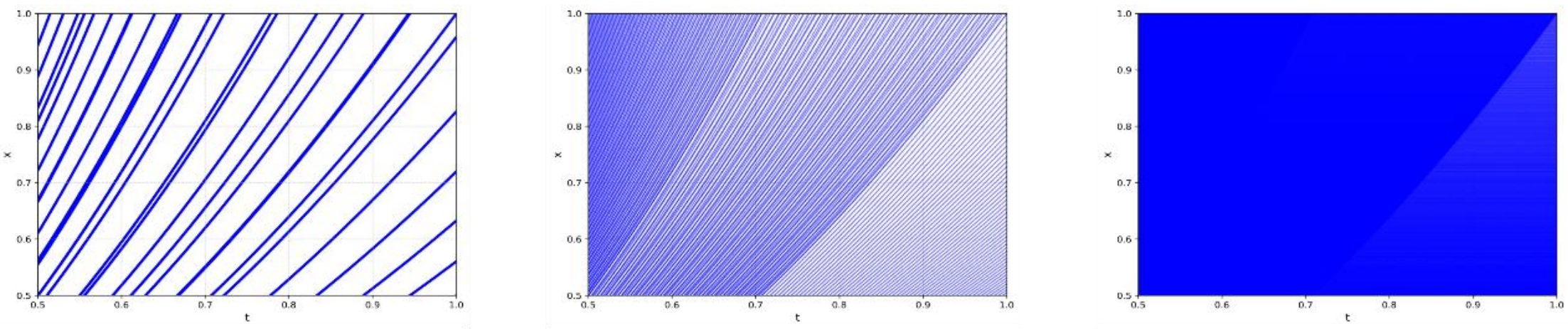

N=10 N=100 N=1000

Fig. 3. Distribution of base points in the solution domain obtained by sampling 10, 100, and 1000 points on the initial and boundary conditions, respectively.

## 3.2 Lie-AS-ISC

The base points obtained from the above procedure are denoted by $\{(x_i, t_i, u_i)\}_{i=1}^{N_l}$ .Since the invariance surface condition holds only for group-invariant solutions, this paper focuses its discussion on group-invariant solutions. Prior to training, a set of collocation points $\{(x_i, t_i)\}_{i=1}^{N_r}$ is sampled within the domain to enforce the PDE residual constraints. After training the neural network, predictions at the base points are obtained, yielding the corresponding approximate solution values. To enhance sampling effectiveness, a dual-criterion selection strategy is employed. The first criterion is the absolute error between the predicted and exact solutions, $|u_{pred} - u_{exact}|$, and the second is the absolute value of the invariant surface condition, $|\tau - \xi u_x - \eta u_t|$. Points with larger combined indicator values are assigned higher selection probabilities, thereby prioritizing regions that are difficult to approximate.

From the base points, $M$ points are selected according to this strategy, denoted by $\left\{\left(x_j,t_j,u_j\left(x_j,t_j\right)\right)\right\}_{j=1}^{M} v$.

These selected points are then incorporated into the loss function:

$$L_{lie}=\frac{1}{M}\sum_{j=1}^{M}\left|\hat{u}\left(x_j,t_j\right)-u_j\left(x_j,t_j\right)\right|^2, \tag{29}$$

$$L_{ISC}=\frac{1}{M}\sum_{j=1}^{M}\left|h\left(x_j,t_j\right)\right|^2, \tag{30}$$

Here, the invariant surface condition is defined as $h\left(x_j,t_j\right):=\tau-\xi u_x-\eta u_t$. The total loss function is formulated as

$$L=\lambda_r L_r+\lambda_i L_i+\lambda_b L_b+\lambda_{lie} L_{lie}+\lambda_{ISC} L_{ISC}, \tag{31}$$

where $\lambda_r,\lambda_i,\lambda_b,\lambda_{lie}$ and $\lambda_{ISC}$ are weighting coefficients corresponding to the residual, initial condition, boundary condition, Lie symmetry, and invariant surface condition losses, respectively. The weighting parameters are adaptively adjusted using a negative log-likelihood-based strategy [36].

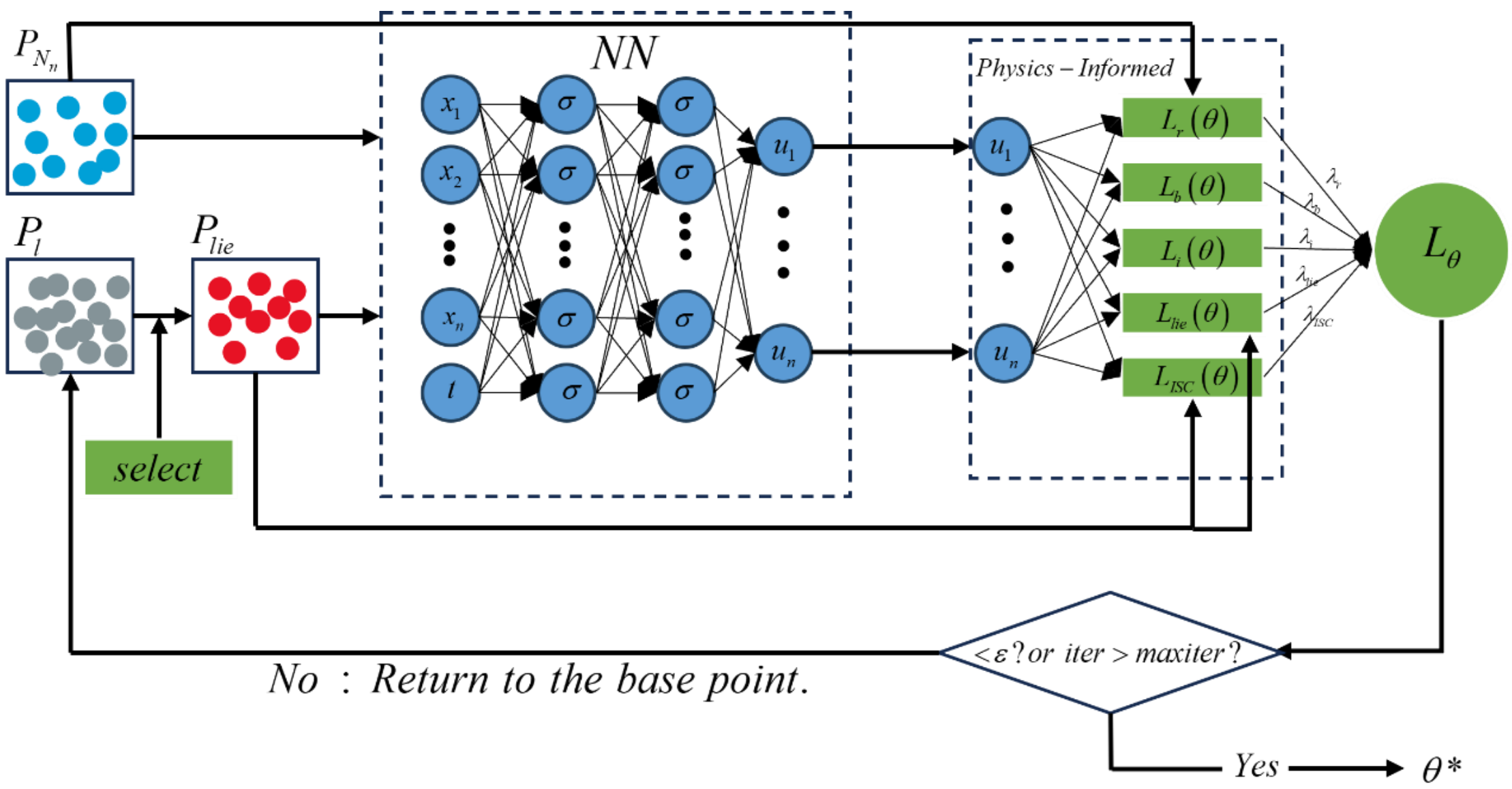

Fig. 4. Overview of the Lie-AS-ISC framework.

Fig. 4 illustrates the framework of the proposed Lie-AS-ISC method, a Lie symmetry-based adaptive sampling physics-informed neural network (PINN). Here, $P_{N_n}$ denotes the randomly sampled collocation points used to enforce the residual loss $L_r$, while $P_l$ represents the global labeled points generated via Lie symmetry transformations. The set $P_{lie}$ consists of selected points from $P_l$, determined based on the Lie symmetry loss $L_{lie}$ and the invariant surface condition loss $L_{ISC}$. In the initial training stage, only $P_{N_n}$ is used, and the loss function includes three components: $L_r, L_b$ and $L_i$. After screening, a subset $P_{lie}$ is obtained and incorporated into subsequent training. The loss function is then extended to include five components: $L_r, L_b, L_i, L_{lie}$ and $L_{ISC}$. The set $P_{lie}$ is updated iteratively after each training epoch. This iterative refinement mechanism enables adaptive sampling in regions with complex solution structures, thereby improving both accuracy and convergence efficiency.

## 4. Numerical results

Lie-S is derived from Lie-AS-ISC by omitting the adaptive sampling component; specifically, it retains only the points obtained via symmetric transformation sampling, without performing subsequent adaptive resampling. Lie-AS is derived from Lie-AS-ISC by removing the invariant surface condition; specifically, the invariant surface condition is eliminated from both the resampling metric and the loss function.

### 4.1 Klein-Gordon equation

We first consider the Klein–Gordon equation in the form:

$$\begin{cases} u_{tt} - u_{xx} + u + u^2 = f(x,t), x \in [0,1], t \in [0,1] \\ u(0,x) = \varphi(x) \\ u(t,0) = \phi^l(t), u(t,1) = \phi^r(t) \end{cases}, \quad (32)$$

The exact solution is given by

$$u(x,t) = 15\cos(\pi(x-3t)) + (x-3t)^2, \quad (33)$$

This equation admits a translational symmetry, with the corresponding transformation:

$$x^* = x + 3\varepsilon, t^* = t + \varepsilon, u^* = u, \tag{34}$$

From this transformation, the relationship among $x^*, t^*$ and $u^*$ can be derived as:

$$\begin{cases} x^* = 3t^* + x - 3t \\ u^* = u \end{cases}, \tag{35}$$

Accordingly, the invariant surface condition is obtained as:

$$3u_x + u_t = 0. \tag{36}$$

Based on the exact solution, the corresponding initial and boundary conditions are derived. The model is trained using a combination of the Adam and L-BFGS optimizers. Specifically, Adam is first applied for 5,000 iterations with a learning rate of 0.001, followed by L-BFGS optimization for 50,000 iterations with a learning rate of 0.5. Resampling is performed three times during training. The neural network consists of two hidden layers, each with 50 neurons. The number of interior collocation points is set to $N_r = 1500$, while $M = 500$ points are selected for adaptive sampling. In addition, 200 boundary points and 100 initial points are used.

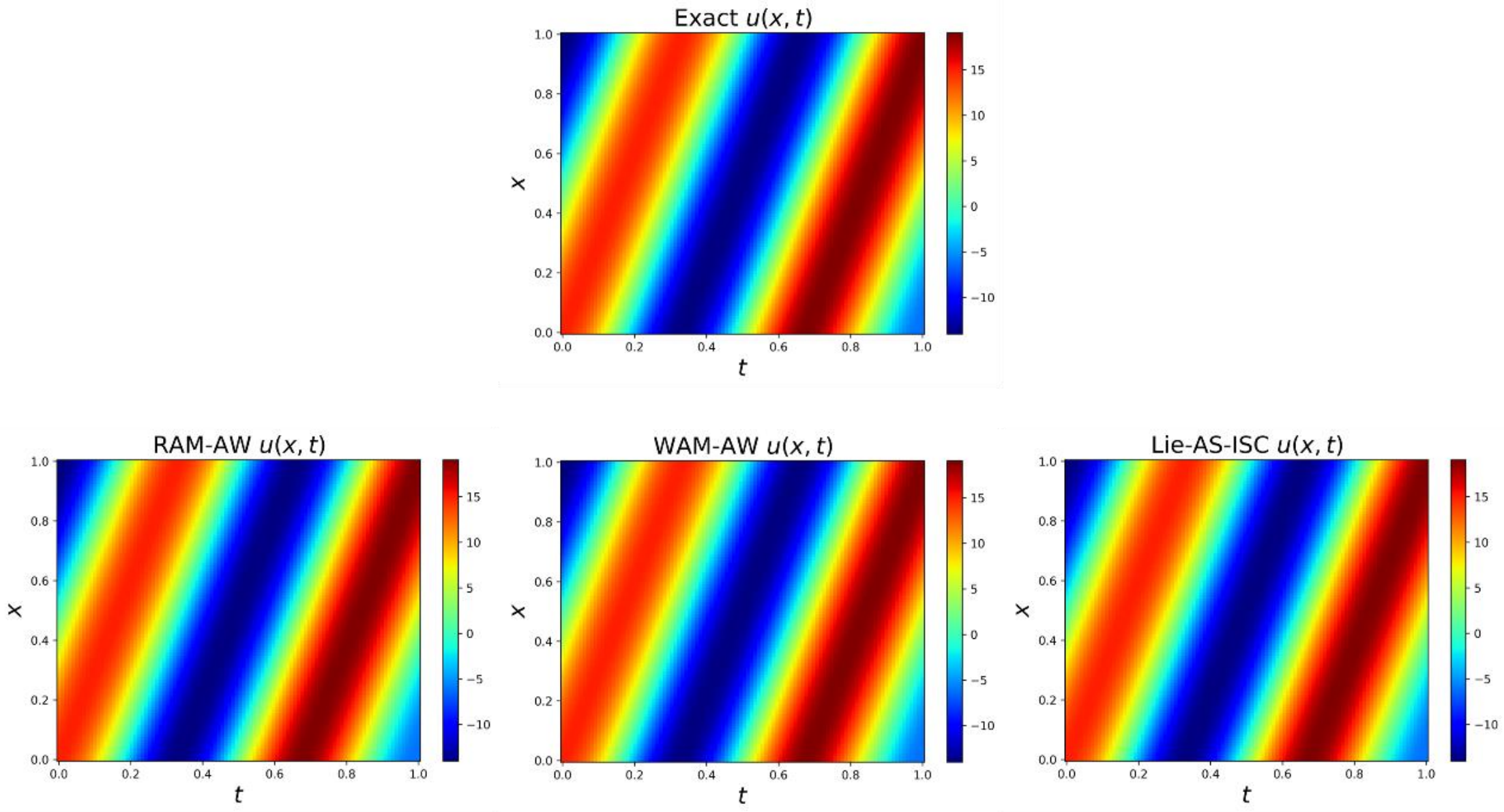

Fig. 5. Comparison of exact and predicted solutions for the Klein–Gordon equation using RAM-AW, WAM-AW, and the proposed Lie-AS-ISC method.

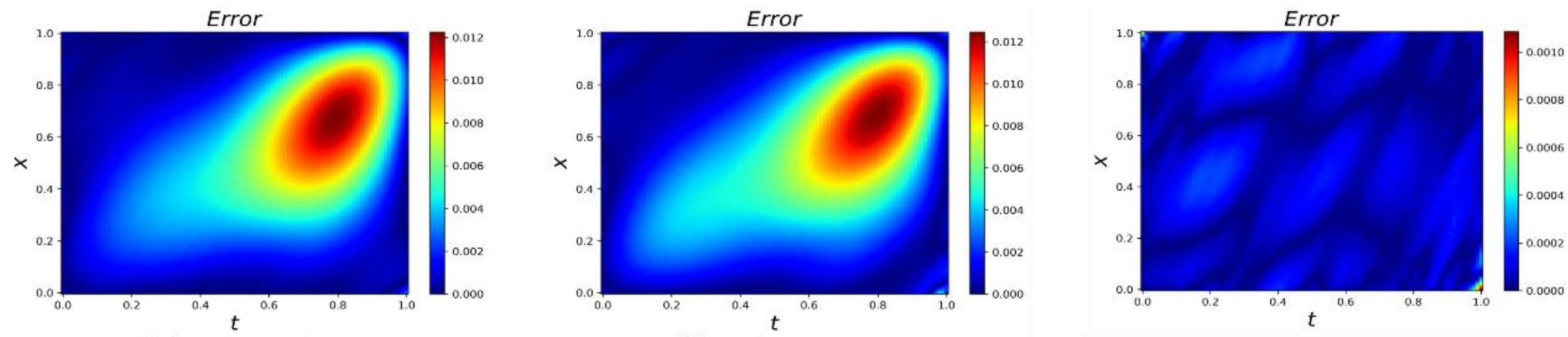

(a)RAM-AW (b)WAM-AW (c)Lie-AS-ISC

Fig. 6. Comparison of absolute errors for the Klein–Gordon equation among RAM-AW, WAM-AW, and the proposed Lie-AS-ISC method.

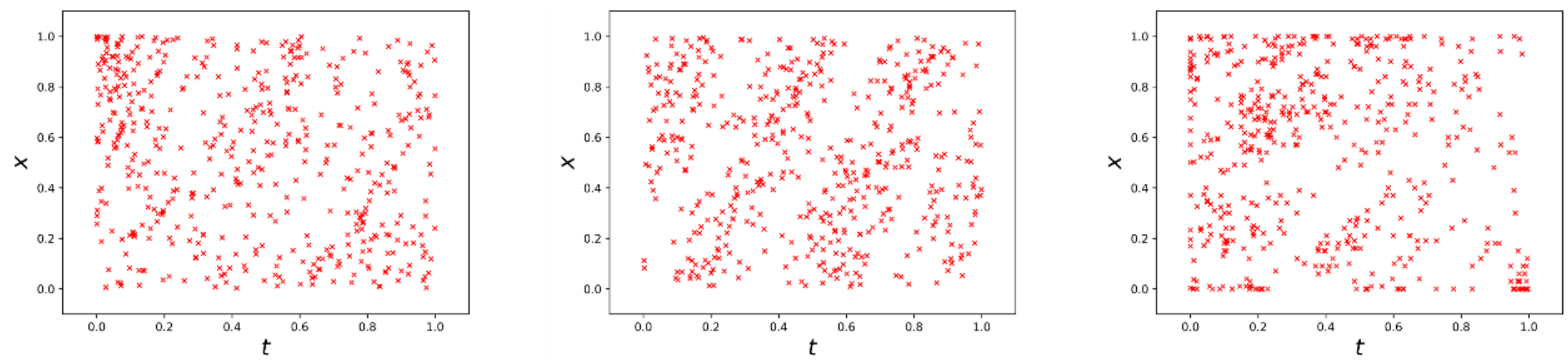

(a)RAM-AW (b)WAM-AW (c)Lie-AS-ISC

Fig. 7. Comparison of resampling point distributions for the Klein–Gordon equation among RAM-AW, WAM-AW, and the proposed Lie-AS-ISC method.

RAM-AW and WAM-AW are residual-based and gradient-based adaptive sampling methods, respectively [17], while APINN is a domain decomposition approach based on gate-controlled networks [37]. Table 1 presents a comparison of the L2 errors between these methods and the proposed approach.

As shown in Fig. 6, both RAM-AW and WAM-AW exhibit significant regional errors, particularly in challenging regions. In contrast, the proposed Lie-AS-ISC method achieves consistently low error across the entire domain, without noticeable regions of large deviation. Quantitative results in Table 1 further demonstrate that, after three rounds of resampling, the error of Lie-AS-ISC is nearly two orders of magnitude lower than that of the competing methods. At the same time, Lie-AS-ISC demonstrates clear advantages over Lie-AS and Lie-S.

Table 1 Klein-Gordon equation: relative L2 error of three round resampling

| Klein-Gordon | 1 | 2 | 3 |
|---|---|---|---|
| RAM-AW | 2.80e-04 | 4.36e-04 | 3.85e-04 |
| WAM-AW | 2.90e-04 | 4.59e-04 | 3.99e-04 |
| APINN | 5.68e-03 | / | / |
| Lie-S | 1.40e-03 | / | / |
| Lie-AS | 1.96e-04 | 5.47e-05 | 4.31e-05 |
| Lie-AS-ISC | 1.82e-05 | 1.22e-05 | 7.52e-06 |

## 4.2 KdV equation

We consider the KdV equation in the form::

$$\begin{cases} u_t + uu_x + u_{xxx} = f(x,t), x \in [0,1], t \in [0,1] \\ u(0,x) = \varphi(x) \\ u(t,0) = \phi^l(t), u(t,1) = \phi^r(t) \end{cases}, \tag{37}$$

The exact solution is given by

$$u(x,t) = (x-2t)^2 + 20\sin(\pi(x-2t)), \tag{38}$$

This equation admits a translational symmetry with the transformation:

$$x^* = x + 2\varepsilon, t^* = t + \varepsilon, u^* = u, \tag{39}$$

Accordingly, the invariant surface condition is:

$$2u_x + u_t = 0, \tag{40}$$

The relationship among $x^*, t^*$ and $u^*$ can be expressed as:

$$\begin{cases} x^* = 2t^* + x - 2t \\ u^* = u \end{cases}, \tag{41}$$

The model is trained using a combination of the Adam and L-BFGS optimizers. Specifically, Adam is applied for 5,000 iterations, followed by L-BFGS for 50,000 iterations. After each training stage, 500 points are resampled from the base points. The neural network consists of two hidden layers, each with 50 neurons. After three rounds of resampling, the error of the proposed Lie-AS-ISC method is reduced to 8.54e-06, which is at least two orders of magnitude lower than that of residual- and gradient-based adaptive sampling methods.

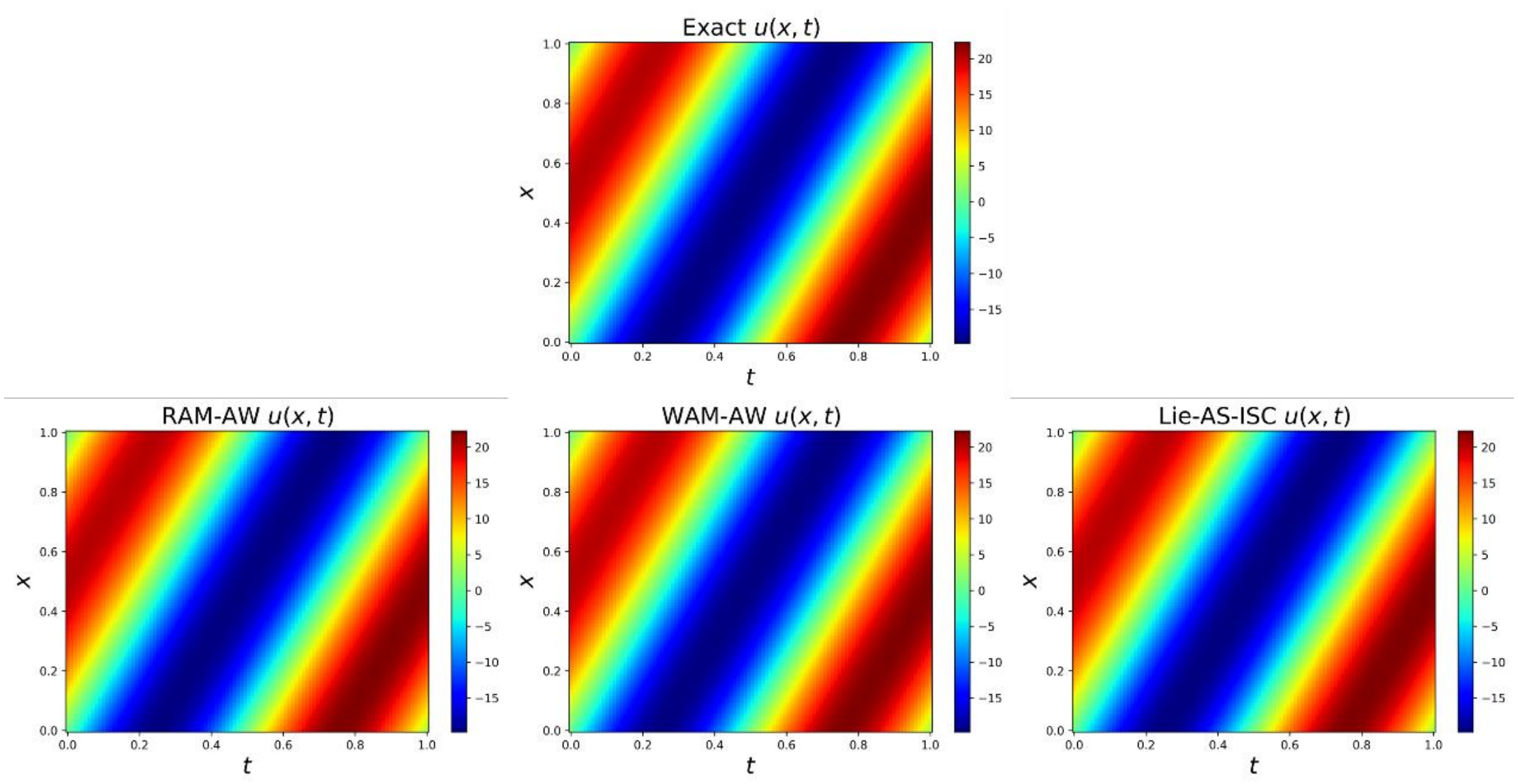

Fig. 8. KDV: exact and predicted solutions

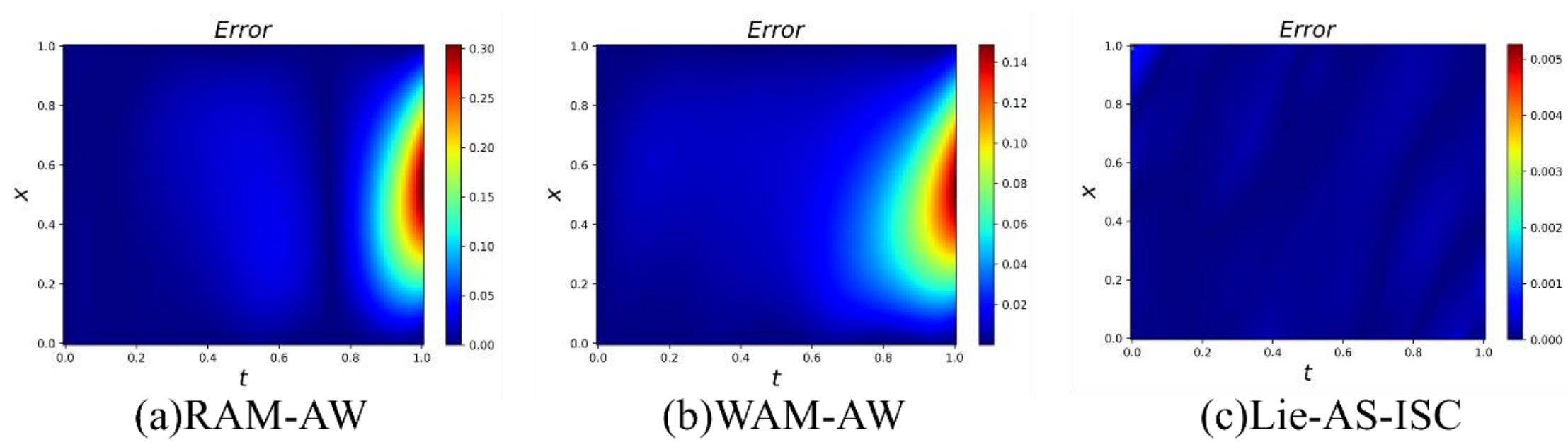

Fig. 9. KDV: Absolute Error

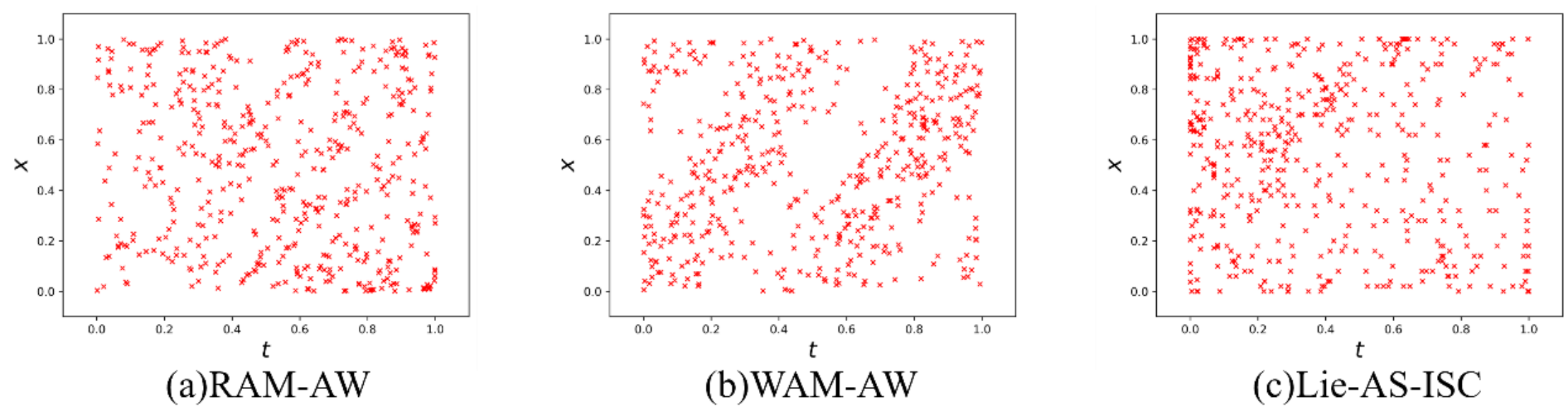

Fig. 10. KDV: Distribution of resampling points

As shown in Fig. 9, both RAM-AW and WAM-AW exhibit significant errors in the right region of the domain. In contrast, the proposed Lie-AS-ISC method achieves lower global error and effectively suppresses local error concentrations. Quantitative results in Table 2 further demonstrate that the L2 error of Lie-AS-ISC is substantially lower than that of the comparison methods, achieving an improvement of nearly three orders of magnitude. Lie-AS-ISC significantly outperforms Lie-AS and Lie-S.

Table 2 KDV equation: Relative L2 error after three rounds of resampling

| KDV | 1 | 2 | 3 |
|---|---|---|---|
| RAM-AW | 9.62e-03 | 3.83e-03 | 3.60e-03 |
| WAM-AW | 3.82e-03 | 2.07e-03 | 2.09e-03 |
| APINN | 2.86e-01 | / | / |
| Lie-S | 1.46e-04 | / | / |
| Lie-AS | 7.43e-05 | 2.53e-05 | 1.69e-05 |
| Lie-AS-ISC | 3.33e-05 | 1.05e-05 | 8.54e-06 |

## 4.3 Inverse Problem of the KdV Equation

Traditional numerical approaches for solving inverse partial differential equation (PDE) problems often suffer from strong dependence on observational data and limited robustness. In contrast, physics-informed neural networks (PINNs) incorporate both observational data and the underlying physical laws into the learning process, providing a flexible and efficient framework for inverse problem solving. To demonstrate the effectiveness of the proposed Lie-AS-ISC method in inverse settings, we consider the inverse problem associated with the KdV equation, which is formulated as follows:

$$\begin{cases} u_t + kuu_x + u_{xxx} = f\left(x,t\right), x \in \left[0,1\right], t \in \left[0,1\right] \\ u\left(0,x\right) = \varphi\left(x\right) \\ u\left(t,0\right) = \phi^l\left(t\right), u\left(t,1\right) = \phi^r\left(t\right) \end{cases},$$

Our goal is to recover the unknown parameter $k$together with the solution of the governing equation.

The training configuration is as follows. A fully connected neural network with two hidden layers, each containing 50 neurons, is employed. A total of 500 observation points are randomly sampled from the computational domain. The network is initially trained using the Adam optimizer for 5,000 iterations, followed by further optimization with the L-BFGS algorithm for 10,000 iterations. In addition, the adaptive resampling procedure is performed three times during training.

The exact value of the parameter k is 2. The proposed Lie-AS-ISC method accurately identifies the parameter with an estimated value of 2.000026, whereas WAM-AW fails to achieve convergence. Furthermore, Lie-AS-ISC attains an $L_2$ error of only 6.75e-06 for the predicted solution, while the corresponding error produced by WAM-AW reaches 9.92e-01. These results demonstrate the superior accuracy and robustness of Lie-AS-ISC in solving the inverse KdV problem.

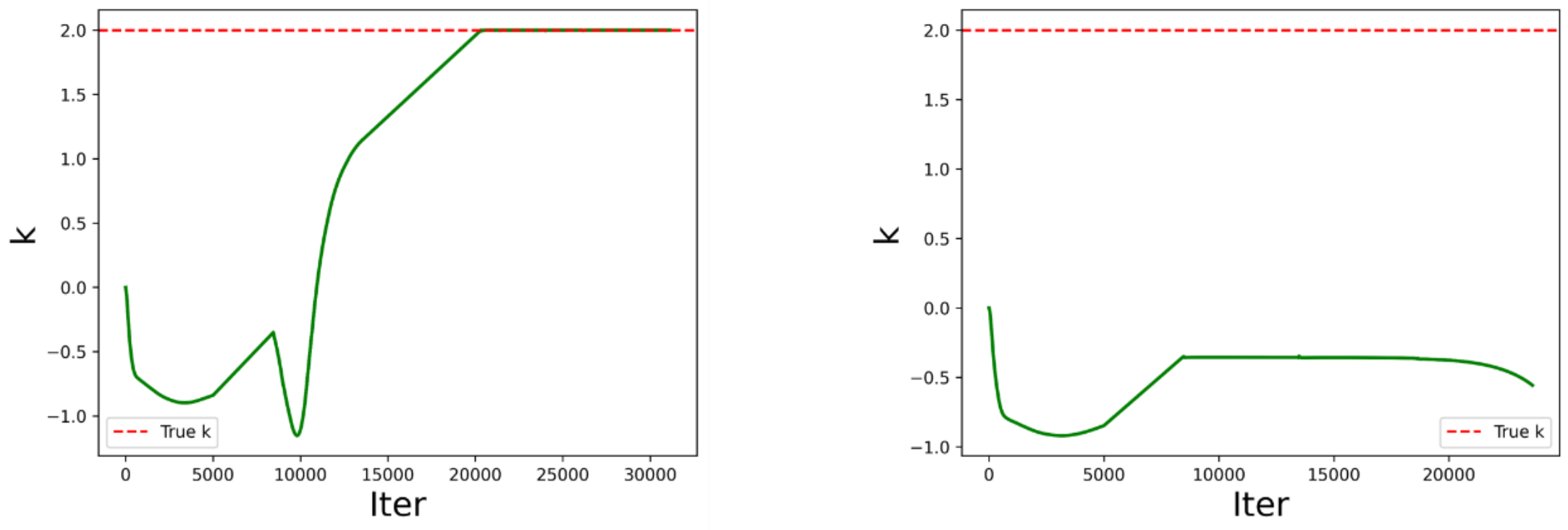

Fig. 11 Left: Convergence history of the parameter $k$ obtained by Lie-AS-ISC; Right: Convergence history of the parameter $k$ obtained by WAM-AW.

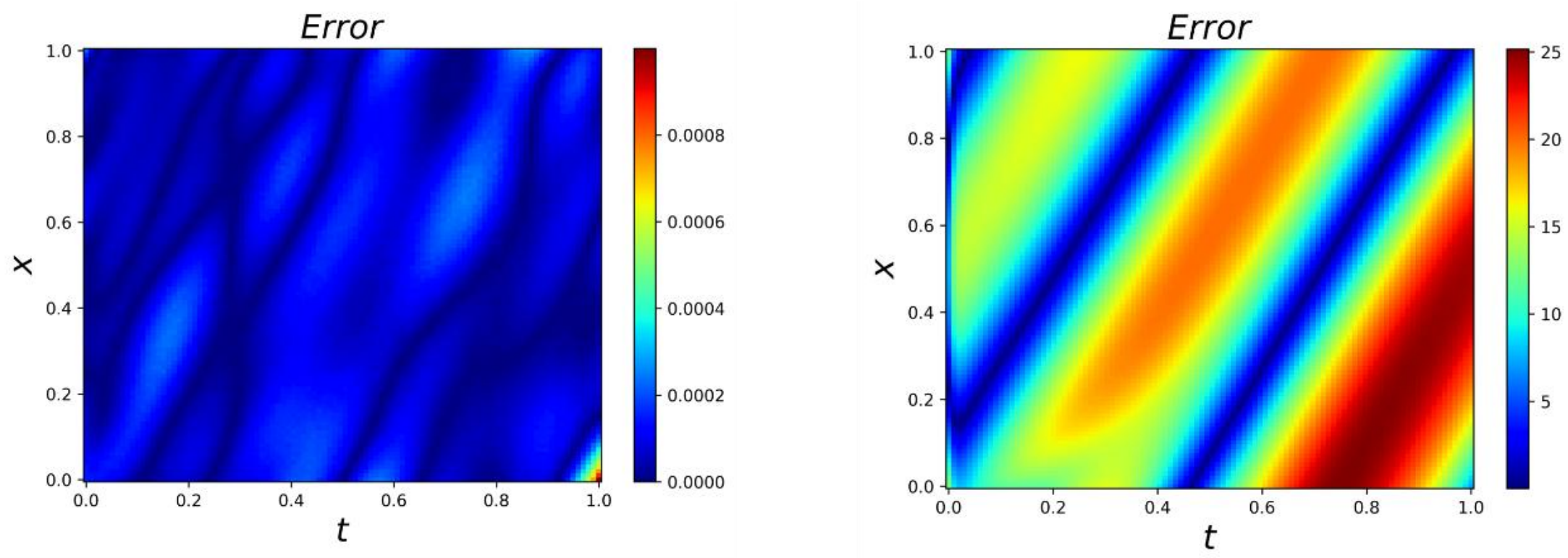

Fig. 12 Left: $L_2$ error of the solution obtained by Lie-AS-ISC; Right: $L_2$ error of the solution obtained by WAM-AW.

## 5. Conclusion

This paper investigates partial differential equations with Lie symmetry groups and proposes a Lie symmetry-based adaptive sampling physics-informed neural network, termed Lie-AS-ISC. First, the exact solution and the functional relationships among variables are derived via Lie symmetry transformations of the governing equations. These relationships are then used to construct function mappings that effectively cover the entire solution domain, thereby generating base points with exact solution values. Subsequently, after each training stage, new sampling points are selected based on a dual-criterion strategy that incorporates both absolute error and invariant surface conditions. These points are further constrained through a corresponding dual-loss mechanism. In addition, to balance multiple loss terms, a negative log-likelihood-based adaptive weighting method is employed to ensure stable and effective training. To evaluate the performance of the proposed method, extensive experiments are conducted across three types of symmetries: translational, rotational, and scaling symmetries. Benchmark problems, including the Klein–Gordon equation and KdV equation, are considered. The results demonstrate that the proposed method significantly reduces prediction errors and effectively mitigates the issue of poor approximation in challenging regions commonly encountered in conventional PINNs. Future work will focus on further integrating intrinsic properties of partial differential equations, such as Lie symmetries, with deep learning frameworks to enhance model robustness and generalization.